\documentclass[11pt]{article}
\usepackage{amssymb}
\usepackage{graphicx}
\parskip \baselineskip
\newcommand{\ds}{\displaystyle}

\newcommand{\Z}{\mathbb{Z}}

\newcommand{\ol}{\overline}

\newcommand{\ra}{\rightarrow}
\newcommand{\Ra}{\Rightarrow}

\date{}

\begin{document}

\title{Output-polynomial enumeration of all fixed-cardinality ideals \\  of a poset, respectively all fixed-cardinality subtrees of a tree}

\author{Marcel Wild}

\maketitle

\begin{abstract}
The $N$ cardinality $k$ ideals of any $w$-element poset ($k \leq w$ fixed) can be enumerated in time $O(Nw^3)$. The corresponding bound for $k$-element subtrees of a $w$-element tree is $O(Nw^5)$. An algorithm is described that by the use of wildcards displays all order ideals of a poset in a compact manner, i.e. not one by one.
\end{abstract}

\section{Introduction}

An {\it implication} is a Boolean formula of type 

(1) \qquad $(a_1 \wedge a_2 \wedge \cdots \wedge a_m) \ra (b_1 \wedge b_2 \wedge \cdots \wedge b_n)$.

Each finite closure system ${\cal C}$ can be viewed as the set of models of a suitable family $\Sigma$ of implications, called an {\it implicational base}. Conversely, each $\sum$ triggers a closure system ${\cal C}$. It was known that from $\Sigma$ one can enumerate ${\cal C}$ in output-polynomial time.  Here we show that for any fixed $k$ the $k$-element members of ${\cal C}$ can be enumerated in output-polynomial time as well, provided ${\cal C}$ satisfies some technical condition (Theorem 2). As Corollary 3 and Corollary 4 we obtain the two results mentioned in the abstract. 

Before we come to the section break up, let us give some background on Corollary 3. For an arbitrary $w$-element poset $(W, \leq)$ we shall write $N$ for the number of its (order) ideals, and $N_k$ for the number of $k$-element ideals $(0\leq k \leq w)$. Calculating $N$ is a \#$P$-complete problem, but for special cases, such as $2$-dimensional posets or interval orders, one can do it in $w$-polynomial time. As to enumerating ($=$ producing) all $N$ ideals, obviously this can't be done in $w$-polynomial time since for general posets (e.g. antichains) $N$ grows exponentially with $w$. However, it can be done in {\it output-polynomial\footnote{Thus in time $O(N^cw^d)$ for some $c,d\in \Z_+$, as opposed to say $O(N^{\sqrt{w}})$ or $O(N^{\log N})$ which are not output-polynomial.}} time. In fact Medina and Nourine [MN] enumerate the ideals in time $O(Nd)$ where $d$ is the maximum number of covers of an element of $W$.

Let us turn to ideals of fixed cardinality $k$. Calculating $N_k$ can be done in $w$-polynomial time for interval orders [S], and for some posets explicite closed formulae for $N_k$ are given in [C]. As to enumerating the $N_k$ ideals of cardinality $k$, for interval orders it was known to be doable in $N_k$-linear time [HNS]. We show in Corollary 3 that $N_k$-linear time can be achieved for general posets.

Sections 2 and 3 introduce the $0,1,2$-algorithm which underlies the proof of Theorem 2 in Section 4. The $0,1,2$-algorithm is convenient for that particular proof, but gets fine-tuned in Section 5 to the $(a,b)$-algorithm which however targets only {\it singleton premise} implications (so $m=1$). This kind of implication, which is intimately linked to posets, is the sole kind we are concerned with in the remainder of the article. Applications to stable marriages, modular lattices, and others are pointed out.

The number $N_k$ of $k$-element ideals of a poset $(W, \leq)$ can be viewed as the $k$-th {\it Whitney number} of the associated distributive lattice $D$. In Section 6 we pit the $(a,b)$-algorithm against a straightforward recursive method to get the Whitney numbers.  Section 7 investigates variations of the $(a,b)$-algorithm, for instance whether it carries over from implications like $a \ra (b_1 \wedge b_2 \wedge b_3)$ to other kinds of ``quasi-implications'' such as $a \ra (\ol{b}_1 \wedge \ol{b}_2 \wedge \ol{b}_3)$.  

\section{The $0,1,2$-algorithm by example}
Rather than (1) we shall mostly adopt the set theoretic notation $\{a_1, \cdots, a_m\} \ra \{b_1, \cdots, b_n \}$ for implications. If both the {\it premise} $A = \{a_1, \ldots, a_m\}$ and {\it conclusion} $B = \{b_1, \ldots, b_n\}$ are subsets of our universe $W$ we say that $X \subseteq W$ {\it satisfies} $A \ra B$ if $A \subseteq X$ implies $B \subseteq X$. In other words, either $A \not\subseteq X$ or $A \cup B \subseteq X$. Thus $X$ corresponds to a {\it model} ($=$ satisfying truth value assignment) of formula (1). In the sequel we may assume that $A \cap B = \emptyset$ since $A \ra B$ is equivalent to $A \ra (B \backslash A)$.

We introduce the $0,1,2$-algorithm on a toy example with $W = [7]:= \{1,2, \ldots, 7\}$. Consider the family of $h=4$ implications
\begin{center}
$\sum \quad : = \quad \{A_1 \ra B_1, A_2 \ra B_2, A_3\ra B_3, A_4 \ra B_4\}$

$ = \quad \{\{5\} \ra \{6,7\}, \quad \{6\} \ra \{3 \}, \quad \{1,2,3\} \ra \{7\}, \quad \{3\} \ra \{4,5\} \}$.
\end{center}
In order to calculate the family Mod $=$ Mod$(\sum)$ of all $\sum$-{\it models} (i.e. those sets $X \subseteq W$ that are models of {\it all} $(A_i \ra B_i) \in \sum$) we represent subsets of $W$ by their characteristic bitstrings of length $w = |W| = 7$ but with the additional proviso of a don't care symbol 2. Thus the powerset ${\cal P}(W)$ is $(2,2, \cdots, 2)$, and in general a $\{0,1,2\}$-{\it valued row} like $r = (0,2,1,2,2,0,2)$ can be viewed as an interval in ${\cal P}(W)$, i.e. $r = \{X \in {\cal P}(W): \ \{3\} \subseteq X \subseteq \{2,3,4,5,7\}\}$.  Coupled to each $\{0,1,2\}$-valued row are three obvious subsets of $W$, which for $r$ above are zeros$(r) = \{1,6\}$, ones$(r) = \{3\}$, twos$(r) = \{2,4,5,7\}$.

Generally we put Mod$_0: = {\cal P}(W)$ and for $1 \leq i \leq h$ let Mod$_i$ be the family of all $X \subseteq W$ that satisfy the first $i$ implication of $\sum$. In particular Mod$_h =$ Mod. Thus if 

(2) \qquad ${\cal P}_i : = \{ X \in {\cal P}(W): \ X \ \mbox{satisfies} \ A_i \ra B_i \},$

then Mod $= {\cal P}_1 \cap {\cal P}_2 \cap \cdots \cap {\cal P}_h$.
The main idea is to calculate Mod$_{i+1}$ from Mod$_i$ by {\it excluding} all $X \in \ \mbox{Mod}_i$ that violate the $(i+1)$-th implication. Such a {\it principle of exclusion} ({\it POE}) is applicable in a wider context, see Section 3.

By ``imposing'' $A_1 \ra B_1$ upon $r_1 = (2,2,2,2,2,2,2)$ in an obvious way we get Mod$_1 = r_1 \uplus r_3$ (Table 1) where $\uplus$ indicates disjoint union.
Rather than imposing $A_2 \ra B_2$ on $r_2$ {\it and} $r_3$ we keep focusing attention on the {\it topmost} row  (here $r_2$) of our {\it working stack}, i.e. we apply the well known {\it last in first out} (LIFO) principle. LIFO entails recording for each row of the working stack which implication is ``pending'' (e.g. the 2nd implication for $r_2$ and the 3rd implication for $r_5$). Imposing $A_2 \ra B_2$ upon $r_2$ again splits that row and results in the working stack (top to bottom) $\{r_4, r_5, r_3\}$. Imposing $A_3 \ra B_3$ upon $r_4$ is more cumbersome because $|A_3| > 1$.  Obviously $r_4 = S^- \uplus S^+$ where
$$S^-: = \{X \in r_4: \ A_3 \not\subseteq X\} \quad \mbox{and} \quad S^+: = \{X \in r_4 : \ A_3 \subseteq X\}.$$
One checks that $S^- = r_6 \uplus r_7 \uplus r_8$ where the disjointness of rows is noteworthy\footnote{Had the first three components of $r_6, r_7, r_8$ been $(0,2,2), (2,0,2), (2,2,0)$ respectively, the union would not have been disjoint.}. All $X \in S^-$ satisfy $A_3 \ra B_3$ but not all $X \in S^+$ satisfy it; for instance $X = (1,1,1,1,0,0,0) \in S^+$ doesn't. However, it is clear that $r_9 
\subseteq S^+$ comprises exactly the sets $X \in S^+$ that satisfy $A_3 \ra B_3$. Imposing $A_4 \ra B_4$ upon $r_6$ results in the row $(0,2, {\bf 0}, 2,0,0,2)$ since $r_6$ has $0$ on its fifth position. All implications having been imposed that row becomes the bottom row of the {\it final stack} (Table 2). Similarly $r_7$ gets its third component 2 switched to $0$ and is put on the final stack. Further $r_8$ is put on the final stack {\it unaltered} but $r_9$ must be {\it deleted} altogether since no $X \in r_9$ satisfies $A_4 \ra B_4$. As to $r_5$, imposing the pending implication $A_3 \ra B_3$ is pointless since again no $X \in r_5$ satisfies $A_4 \ra B_4$, and so $r_5$ must be deleted. The working stack now consists of $r_3$ alone. Imposing the three pending implications on $r_3$ results in $(2,2,1,1,1,1,1)$ which becomes the last row of our final stack. It follows that Mod is the disjoint union of the four $\{0,1,2\}$-valued rows in Table 2, in particular $|\mbox{Mod}| = 4+4+4+8 = 20$.

\newpage

\begin{tabular}{c|c|c|c|c|c|c|c|l} 
& 1 & 2& 3& 4 & 5& 6 & 7 & \\ \hline
& & & & & & & & \\ \hline
$r_1=$ & 2 & 2 & 2 & 2 & 2 & 2 &2 & p. imp 1\\ \hline
& & & & & & & & \\ \hline
$r_2=$ & 2 & 2 & 2 & 2 & ${\bf 0}$ & 2 & 2 & p. imp 2\\ \hline
$r_3=$ & 2 & 2 & 2 & 2 & ${\bf 1}$ & 1 & 1 & p. imp 2 \\ \hline
& & & & & & & & \\ \hline
$r_4 =$ & 2 & 2 & 2 & 2 & 0 & ${\bf 0}$ & 2  & p. imp 3 \\ \hline
$r_5=$ & 2 & 2 & 1 & 2 & 0 & ${\bf 1}$ & 2 & p. imp 3 \\ \hline
$r_3=$ & 2 & 2 & 2 & 2 & 1 & 1 & 1 & p. imp 2 \\ \hline
& & & & & & & & \\ \hline
$r_6=$ & ${\bf 0}$ & ${\bf 2}$ & ${\bf 2}$ & 2 & 0 & 0 & 2 & p. imp 4 \\ \hline
$r_7=$ & ${\bf 1}$ & ${\bf 0}$ & ${\bf 2}$ & 2 & 0  & 0 & 2 & p. imp 4\\ \hline
$r_8=$ & ${\bf 1}$ & ${\bf 1}$ & ${\bf 0}$ & 2 & 0 & 0 & 2 & p. imp 4 \\ \hline
$r_9=$ & ${\bf 1}$ & ${\bf 1}$ & ${\bf 1}$ & 2 & 0 & 0  & 1 & p. imp 4 \\ \hline
$r_5=$ & 2 & 2 & 1 & 2 & 0  & 1 & 2 & p. imp 3 \\ \hline
$r_3=$ & 2 & 2 & 2 & 2 & 1 & 1 & 1 & p. imp 2 \\ \hline
& & & & & & & & \\ \hline
$r_3=$ & 2 & 2 & 2 & 2 & 1 & 1 & 1 & p. imp 2 \\ \hline 
\end{tabular} \hspace*{2cm}
\begin{tabular}{|c|c|c|c|c|c|c|c} 
2 & 2 & 1 & 1 & 1 & 1 & 1 & $\ra 4$\\ \hline
1 & 1 & 0 & 2 & 0 & 0 & 2 & $\ra 4$\\ \hline
1 & 0 & 0 & 2 & 0 & 0 & 2 & $\ra 4$ \\ \hline
0 & 2 & 0 & 2 & 0 & 0 & 2 & $\ra 8$ \\ \hline
\end{tabular}

\hspace*{2cm} Table 1 \hspace*{8cm} Table 2

\section{Formal definition of the $0,1,2$-algorithm and the avoidance of row deletions}

In each version of mentioned POE (more details below) one must show how a general row $r$ that satisfies the first $i$ constraints (thus $r \subseteq \ \mbox{Mod}_i$) either gets deleted or gives rise, upon imposing the $(i+1)$-th constraint, to {\it proper sons} $r'_1, \cdots, r'_s$ such that
$$r \cap \ \mbox{Mod}_{i+1} \quad = \quad r'_1 \uplus r'_2 \uplus \cdots \uplus r'_s.$$
Then LIFO can be applied as illustrated in Section 2 and yields Mod as a disjoint union of rows (final stack). 
In the present article the constraints are implications and in order to impose $A_{i+1} \ra B_{i+1}$ on $r$ we make a case distinction.

Case 1: For all $X \in r$ the disjunction $(A_{i+1} \not\subseteq X$ or $B _{i+1} \subseteq X$) is true. Then $r$ carries over, i.e. $s=1$ and $r'_1 : = r$.

Case 2: The above disjunction fails for all $X \in r$. Then $r$ gets deleted.

The remaining possibility is that $(A_{i+1} \not\subseteq X$ or $B_{i+1} \subseteq X$) holds for some and fails for some $X \in r$. We proceed according to the cardinality of the premise $A_{i+1}$.

Case 3: $|A_{i+1}| = 1$, say $A_{i+1} = \{j\}$. Then the $j$-th entry $\alpha_j$ of $r$ cannot be $\alpha_j=0$ and not all components indexed by $B_{i+1}$ can be 1 (otherwise we'd be in Case 1).

Case 3.1: Some components of $r$ indexed by $B_{i+1}$ are $0$. Then $\alpha_j=2$ (if $\alpha_j$ was $1$, we'd be in Case 2), and so $r'_1  =r_s'$ arises from $r$ by switching $\alpha_j$ to $0$. 

Case 3.2: All components of $r$ indexed by $B_{i+1}$ are either $1$'s or ``anxious'' $2$'s (as if fearing for their freedom).  Subcase 1: $\alpha_j=2$. Then $r'_1$ is obtained by switching $\alpha_j$ to $0$. And $r'_2 = r'_s$ is obtained from $r$ by switching $\alpha_j$ to $1$, along with the anxious $2$'s. Subcase 2: $\alpha_j=1$. Then only one row $r'_1=r'_s$ arises, obtained from $r$ by switching the anxious $2$'s to 1.

Case 4: $|A_{i+1}| > 1$. Analogous to Case 3 all components of $r$ indexed by $A_{i+1}$ are 1 or 2 and not all components  indexed by $B_{i+1}$ are 1. 

Case 4.1: Some components of $r$ indexed by $B_{i+1}$ are $0$. Then not all components or $r$ indexed by $A_{i+1}$ are 1 (otherwise we'd be in Case 2). To fix ideas, assume without loss of generality that $r$ is as in Table 3 and $A_{i+1} \ra B_{i+1}$ is $\{1,2,3,4,5,6\} \ra \{8,9\}$. Thus imagine the first two symbols $\ast$ in $r$ are 1, followed by four 2's. Then Mod$_{i+1} \cap r = \{X \in r: \ A_{i+1} \not\subseteq X\}$ and the latter is $r'_1 \uplus r'_2 \uplus r'_3 \uplus r'_4$, so $s=4$. (Of course the symbols $\ast$ in every $r'_i$ have the same value as the corresponding $\ast$ in $r$.) It is handy to call the boldface $4 \times 4$ pattern in Table 3, or any $t \times t$ pattern of this type, the {\it Flag of Papua} (New Guinea). As in the example of Section 3 where $t =3$, the purpose of the Flag of Papua is to keep the $t$ concerned rows {\it disjoint}.

Case 4.2: All components of $r$ indexed by $B_{i+1}$ are 1 or anxious 2's .

Subcase 1: Not all components of $r$ indexed by $A_{i+1}$ are 1. Say again the first two symbols $\ast$ in $r$ are 1, followed by four $2$'s. However, now Mod$_{i+1} \cap r = S^- \uplus S^+$ with $S^- : = \{X \in r: \ A_{i+1} \not\subseteq X\}$ and $S^+: = \{X \in r: B_{i+1} \subseteq X\}$. As before $S^- = r'_1 \uplus r'_2 \uplus r'_3 \uplus r'_4$ but $S^+ = r'_5$. Subcase 2: All components or $r$ indexed by $A_{i+1}$ are 1. Then Mod$_{i+1} \cap r$ reduces to $S^+ = r_5$. That concludes the formal description of the $0,1,2$-algorithm.

\begin{center}
 \begin{tabular}{l|c|c|c|c|c|c|c|c|c|}
&  1 & 2 & 3 & 4 & 5 & 6 & 7 & 8 & 9 \\ \hline
&    &   &   &   &   &   &   &   & \\ \hline
$r=$ & $\ast$ & $\ast$ & 2 & 2 & 2& 2& $\ast$ & $\ast$ & $\ast$\\ \hline
&    &   &   &   &   &   &   &   & \\ \hline
$r_1'=$ & $\ast$ & $\ast$ & ${\bf 0}$ & ${\bf 2}$ & ${\bf 2}$ & ${\bf 2}$ & $\ast$ & $\ast$ & $\ast$ \\ \hline 
$r_2'=$ & $\ast$ & $\ast$ & ${\bf 1}$ & ${\bf 0}$ & ${\bf 2}$ & ${\bf 2}$ & $\ast$ & $\ast$ & $\ast$\\ \hline 
$r_3'=$ & $\ast$ & $\ast$ & ${\bf 1}$ & ${\bf 1}$ & ${\bf 0}$ & ${\bf 2}$ & $\ast$ & $\ast$ & $\ast$ \\ \hline 
$r_4'=$ & $\ast$ & $\ast$ & ${\bf 1}$ & ${\bf 1}$ & ${\bf 1}$ & ${\bf 0}$ & $\ast$ & $\ast$ & $\ast$ \\ \hline 
$r_5'=$ & 1 & 1 & 1 & 1 & 1 & 1 & $\ast$ & 1 & 1 \\ \hline \end{tabular}

 Table 3
\end{center}



 However, the deletion\footnote{When speaking of row deletions we always mean ``wasteful'' deletions (such as $r_9$ in Sec.2), as opposed to the ``harmless'' deletions where rows give way to proper sons (such as $r_1$ to $r_2, r_3$ in Sec.2).} of rows (Case 2) precludes a theoretic assessment of it.
Fortunately Theorem 1 in [W2] tells us how an ``old'' version of the POE can be purged from row deletions and thus be evaluated.  In brief, it works whenever the existence of a model in a row can be decided fast. Furthermore [W2, Thm.1] also handles the restriction to models of fixed cardinality.  We take the opportunity to sharpen [W2, Thm.1] in the form of Lemma 1 below, but need a few preliminaries.

Some versions of the POE employ additional symbols besides $0,1,2$, say $e, n$, or $a, b$ in Section 5 of the present paper. We then speak of {\it multivalued} rows.  Generalizing ${\cal P}_i$ in (2) other  set families ${\cal P}_i \subseteq {\cal P}(W) \ (1 \leq i \leq h)$, referred to as {\it constraints}, may be concerned. A multivalued row $r$ is {\it feasible} if it contains a {\it model}, i.e. $r \cap {\cal P}_1 \cap \cdots \cap {\cal P}_h \neq \emptyset$.  If  $k \in [w]$ is fixed, then $r$ is {\it extra feasible} if $r$ contains a $k$-element model. The imposition of a constraint on $r$ works as follows. Row $r$ splits into at most $s$ {\it candidate sons}. some of which get killed or altered and the remaining rows are the {\it proper sons}. Getting the cadidate sons is easier than the proper sons. For instance, if the $0,1,2$-algorithm is run on a family of implications $\sum$ then

(3) \qquad $s = \max \{|A| +1 : \ (A\ra B) \in \sum \}$.

{\bf Lemma 1:} Let $W$ be a set of cardinality $w$ and let ${\cal P}_i \subseteq 2^W$ be $h$ constraints.
Fix $k \in [w]$. Suppose some ``old'' version of POE can be employed to produce
disjoint multivalued rows whose elements are the $N > 0$ $k$-element models\footnote{Here and in Sec.5 we prefer $N$ over $N_k$ to denote the number of $k$-element models. If $N=0$ (whence $R =0$) the expression $O(\cdots)$ in Lemma 1 becomes $0$, despite the fact that it takes time to discover there are no $k$-element models. We could have covered the $N=0$ case by writing $R+1$ for $R$ and $N+1$ for $N$ within $O(\cdots)$ but that got too clumsy, particularly in the proof. Note that $N+1$ fits better in Theorem 2, and in Corollary 3 and 4 one always has $N \neq 0$.}. Further assume that for the functions $f(h,w)$ 
and $f^\ast(k,w)$ the following holds: 
\begin{enumerate}
\item[(a)]  For each multivalued row $r$ it costs $O(f(h,w))$ to decide whether it is extra feasible.
\item[(b)]  For each multivalued row $r$ it costs $O(\mbox{Card}(r,k)f^\ast (k,w))$ to list (in ordinary
set notation) the sets $X \in r$ with $|X| =k$. 
\end{enumerate}
Then the old version can be adapted to a new one that avoids row deletions and takes time $O(Rhs(w+f(h,w)))$ to deliver the set of $k$-element models as a disjoint union of $R$ multivalued rows. If the $k$-element models must be enumerated one by one (in ordinary set notation) the cost is $O(Nf^\ast (k,w) + Rhs(w+f(h,w))$). (In [W2] $f^\ast (k,w)$ does not appear since it is unnecessarily subsumed by $f(h,w)$.)

{\it Proof.} Suppose the $O(Rhs(w+f(h,w)))$ cost is established. Then it is easy to see that enumerating the $k$-element models one by one costs an additional amount of $O(Nf^\ast (k,w))$. Namely, it follows at once from (b) in view of the fact that the sum of the $R$ numbers Card$(r,k)$, when $r$ ranges over the $R$ final rows, is $N$.

Showing the $O(Rhs(w+f(h,w)))$ claim requires us to browse the proof of [W2, Theorem 1]. Its overall structure remaining the same, we merely identify two crisp spots where one can save on time.
First, in (13) of [W2] it is postulated that generating the candidate sons of a length $w$ row costs $O(w^2)$. However, in all instances of POE so far one has $s \leq w$ and hence can do with $O(sw)$. Correspondingly the cost $O(w^2)+sO(f(h,w))$ of producing the proper sons of a row (second last line in the proof of [W2, Thm.1]) becomes $O(sw) + sO(f(h,w)) = O(s(w+f(h,w))$. Second, in [W2] it is assumed, merely for cosmetic reasons, that $f(h,w)$ be ``at least linear in $w$.'' The effect is that the cost of producing the proper sons simplifies to $O(wf(h,w))$. As argued in [W2] there are at most $Rh$ occasions where the proper sons of a row are produced. Since no extra costs accumulate from row deletions, the overall cost of producing the $R$ final rows amounts to optionally $O(Rh \cdot wf(h,w))$ (as in [W2]), or to $O(Rh \cdot s(w+f(h,w)))$ which we here prefer.

\section{The main results}

For any family $\sum$ of implications Mod$(\sum) \subseteq {\cal P}(W)$ is well known to be a {\it closure system}, i.e. $W \in \ \mbox{Mod}(\sum$ and from $X, Y \in \ \mbox{Mod}(\sum)$ follows $X \cap Y \in \ \mbox{Mod}(\sum)$. Conversely, given any closure system ${\cal C} \subseteq {\cal P}(W)$, there always are {\it implicational bases} $\sum$ in the sense that ${\cal C}= \mbox{Mod}(\sum)$. For any closure system ${\cal C} \subseteq {\cal P}(W)$ and any set $T \subseteq W$ we write $cl (T)$ for the {\it closure} of $T$, i.e. the smallest $X \in {\cal C}$ with $T \subseteq X$. The {\it length} of an implication $A \ra B$ is $|A| +|B|$, and for a family $\sum$ of implications $||\sum||$ is the sum of the lengths of its members. Observe that $||\sum|| \leq w|\sum|$. If ${\cal C}$ comes along with an implicational base $\sum$, then by [RM, Thm.10.3] there is a $O(||\sum||+w)$ algorithm for computing $cl(T)$ for any $T \subseteq W$.

For a family ${\cal F}$ of closure systems ${\cal C}$ with universes of unbounded cardinality $w$ consider the following {\it disjoint extension property}:

(4) \qquad There is $c\geq 1$ such that for all ${\cal C} \in {\cal F}$, all $Z_0 \in {\cal C}$, all $Y \subseteq W$ and all $k \leq w$, it can be \\
\hspace*{1.3cm}  checked in time $O(w^c)$ whether there is $Z \in {\cal C}$ with $Z_0 \subseteq Z$ and $|Z| =k$ and
$Y \cap Z = \emptyset$.

(If already $Y \cap Z_0 \neq \emptyset$ or if $|Z_0|> k$, the answer is easy.)

We mention that {\it one} possibility for ${\cal F}$ to enjoy (4) is the following. Let ${\cal F}$ consist of convex geometries ${\cal C}$ in the sense of [EJ] which additionally have this property: For every $Z_0 \in {\cal C}$ and every $Y \subseteq W$ with $Y \cap Z_0 = \emptyset$ the set $\{Z \in {\cal C}: Z \supseteq Z_0, Y\cap  Z=\emptyset \}$ has a {\it largest} member $\ol{Z}_0$, and it can be found in time $O(w^c)$. The families ${\cal F}$ in Corollary 3 and 4 are of this kind.

\begin{tabular}{|l|} \hline \\
{\bf Theorem 2:} Let ${\cal F}$ be a family of closure systems that enjoys the disjoint extension \\
property. Then there is an algorithm ${\cal A}$ which for any ${\cal C} \in {\cal F}$, ${\cal C} \subseteq {\cal P}([w])$, that comes  \\
equiped with an $h$-element implicational base $\sum$, and any $k \in [w]$, enumerates the $N$   \\
many $X \in {\cal C}$ with $|X| =k$ in output polynomial time $O((N+1)hs(hw+w^c))$. \\ \\ \hline \end{tabular}

Here $s$ is as in (3). As to $N+1$, see the footnote accompanying Lemma 1.
 
 {\it Proof.} In order to use an enhanced $0,1,2$-algorithm ${\cal A}$ for enumerating all $X \in {\cal C}$ with $|X|=k$ we first verify condition (a) in Lemma 1 for $f(h,w) : =hw+w^c$. So let $r$ be a $\{0,1,2\}$-valued row of length $w$. Compute $Z_0 : = cl(\mbox{ones}(r))$ in time $O(||\sum||+w) = O(hw+w) =O(hw)$. Obviously $Z_0 \subseteq X$ for each $\sum$-model $X \in r$. Hence, if $|Z_0| > k$, then $r$ is not extra feasible. If $|Z_0| \leq k$ put $Y: = \ \mbox{zeros}(r)$. By the disjoint extension property we can test in time $O(w^c)$ whether or not there is a $k$-element $\sum$-model $Z$ with $Z \supseteq Z_0$ and $Z \cap Y = \emptyset$. If yes, then $r$ is extra feasible since $Z \cap Y = \emptyset$ implies $Z \in r$. Vice versa, if the answer is no, then $r$ is not extra feasible. The cost of this test is $O(hw) + O(w^c)$. 
 
 As to Lemma 1(b), putting $\beta = |\mbox{ones}(r)|$ enumerating the $k$-element sets $X \in r$ amounts to enumerate the $(k - \beta)$-element members of the powerset ${\cal P}(\mbox{twos}(r))$. This can be done in time $O(\mbox{Card}(r,k))$, see e.g. Exercise 8 in [K, p.26]. Substituting $f(h,w) = hw+w^c$ and $f^\ast (k,w) =1$ in Lemma 1 yields $O(N+Rhs(hw+w^c)) = O(Nhs(hw+w^c))$. \hfill $\square$
 
\begin{tabular}{|l|}\hline \\
{\bf Corollary 3:} There is an algorithm ${\cal A}$ such that for any $w$-element poset $(W, \leq)$ and any\\
$k \in [w]$ its $N$ many $k$-element order ideals get enumerated by ${\cal A}$ in output polynomial time \\
$O(Nw^3)$.\\ \\ \hline \end{tabular}

{\it Proof.} Associated to each poset $(W, \leq)$ is the closure system ${\cal C}(W, \leq)$ of all its ideals.
We claim that the family ${\cal F}$ of all closure systems ${\cal C}(W, \leq)$ enjoys the disjoint extension property with $c =2$. So let $Z_0$ be an ideal of $(W, \leq)$ and let $Y \subseteq W$ be disjoint from $Z_0$. We may assume that $|Z_0| \leq k$. If $Y' \supseteq Y$ is the (order) filter generated by $Y$ (which remains disjoint from $Z_0$) then each ideal disjoint from $Y$ is contained in the ideal $\ol{Z}_0: = W \backslash Y'$. The set $P : = \ol{Z}_0 \backslash Z_0$ is convex in the sense that $(a< b< c$ and $a,c\in P$) implies $b \in P$. If $|Z_0| +|P| < k$, there is no $k$-element ideal $Z$ extending $Z_0$. On the other hand, if $|Z_0| +|P| \geq k$, we shell any $k - |Z_0|$ elements of $P$ from below. Because $P$ is convex adding these $k-|Z_0|$ elements to $Z_0$ yields a $k$-element ideal $Z$. Calculating $|P|$ costs $O(w^2)$, and so $c=2$.

For any poset $(W, \leq)$ a natural implicational base $\sum = \sum (W, \leq)$ of ${\cal C}(W, \leq)$ is provided by the implications $\{p\} \ra LC(p)$ where $p$ ranges over $W$ and $LC(p)$ is the set of lower covers of $p$. For convenience we admit the minimal elements $p \in W$, albeit they yield trivial implications $\{p\} \ra \emptyset$. Calculating $\sum$ from $(W, \leq)$, say from a $w \times w$ incidence matrix of $\leq$, costs $O(w^2)$. Applying Theorem 2 and noting that $h = |\sum| = w$ and $s=2$, the overall cost of enumerating the $k$-element ideals amounts to $O(w^2) + O(Nhs(hw+w^2))=O(Nw(w^2+w^2))=O(Nw^3)$. \hfill $\square$

 It is an open question whether ``$k$-element ideals'' in Corollary 3 can be generalized to ``ideals of weight $k$'' with respect to some weight function $W \ra \Z_+$. If the answer is affirmative, that would entail (by considering the factor poset of strong components) that for any family $\Sigma$ of singleton-premise implications (possibly with directed cycles) the $\Sigma$-closed $k$-element sets can be enumerated in output-polynomial time. We mention that finding {\it one} maximum weight $k$-element ideal of a poset is NP-hard [FK].

 \begin{tabular}{|l|} \hline \\
 {\bf Corollary 4:} There is an algorithm ${\cal A}$ such that for any $w$-element tree $T$ and any
 $k \in [w]$\\ 
 the $N$
  many $k$-element subtrees of $T$ get enumerated by ${\cal A}$ in time 
  $O(Nw^5)$. \\ \\ \hline \end{tabular}
 
 {\it Proof.} For convenience we identify $T$ and all occuring (induced) subforests of $T$ with their underlying vertex sets. Associated to each tree $T$ is the closure system ${\cal C}(T)$ of all its subtrees. In order to show that the family ${\cal F}$ of all ${\cal C}(T)$ enjoys the disjoint extension property with $c=1$, let $Z_0 \subseteq T$ be a subtree and let $Y \subseteq T$ be a set disjoint from $Z_0$. Then $F = V \backslash Y$ is a subforest, and there is a $k$-element subtree extending $Z_0$ if and only if the connected component $\ol{Z}_0$ of $F$ containing $Z_0$ has cardinality $\geq k$. Finding the size of  the connected component of a vertex in a graph can be done in linear time with respect to the number of edges, hence in our case $O(w)$. This proves $c=1$.
 
 For any tree $T$ an implicational base $\sum$ of ${\cal C}(T)$ is provided by the implications $\{a_1, a_2\} \ra \{b_1, \cdots, b_n\}$ where $a_1, a_2$ are non-adjacent vertices of $T$ and $a_1, b_1, \cdots, b_n, a_2 \ (n \geq 1)$ is the unique path from $a_1$ to $a_2$. One can construct $\sum$ in time $O(w^2)$. Applying Theorem 2 and noting that $h = O(w^2)$ and $s=3$, the overall cost of enumerating the $k$-element subtrees is $O(w^2)+ O(Nhs(hw+w))=O(Nw^2(w^2w+w)) = O(Nw^5)$. \hfill $\square$

There is some literature on counting or enumerating various types of {\it binary} trees and their subtrees. As for arbitrary trees $T$, {\it all} subtrees of $T$ can be enumerated in output-polynomial time [R] but Corollary 4 seems to be new. Observe that different from posets a maximum weight $k$-subtree of $T$ can be found [FK] in time $O(w^4)$.

\section{The $(a,b)$-algorithm} 
 
 In Section 2 imposing the implication $\{1,2,3\} \ra \{7\}$ upon $r_4 = (2,2,2,2, 0,0,2)$ resulted in the Flag of Papua made up by $r_6, r_7, r_8$ and the ``$1$-filler'' row $(1,1,1,2,0,0,1)$. It seems more economic to replace $r_6, r_7, r_8$ by the single row $(n, n, n, 2, 0,0,2)$, where by definition the wildcard $nnn$ means ``at least one $0$ here''.  The resulting {\it implication $n$-algorithm} [W2] indeed beats the $0,1,2$-algorithm in practise but for the proof of Theorem 2 the $0,1,2$-algorithm was more convenient. 
 
For singleton-premise implications the dichotomy $(n, \cdots, n ) \leftrightarrow (1, \cdots, 1)$ boils down to $0 \leftrightarrow 1$ (hence nothing new), but there is another way to improve upon the $0,1,2$-algorithm in this case. In brief,
we shall employ two new symbols $a,b$. Rather than imposing say $\{5\} \ra \{6,7\}$ on $(2,2,2,2,2, 2, 2)$ by splitting it in $r_2, r_3$ (Table 1), we replace $r_1$ by  a single row  $(2,2,2,2,a, b, b)$ 
which {\it by definition} represents the intended set system. 

  Before we embark on the ensuing algorithm, we formally define an {\it $\{0,1,2,a,b\}$-valued} row $r$ 
  as a partition of $[w]$ into four (possibly empty) parts zeros$(r)$, ones$(r)$, twos$(r)$, implications$(r)$, such that if implications$(r) \neq \emptyset$, it is further partitioned as 
$$\mbox{implications}(r)\quad = \quad \mbox{prem}[1] \cup \ \mbox{conc}[1]\cup \cdots \cup \ \mbox{prem}[t] \cup \ \mbox{conc}[t]\quad (t\geq 1),$$
where all prem$[i]$ are singletons, and all conc$[i]$ are nonvoid $(1 \leq i \leq t)$. A set $X \subseteq [w]$ by defintion belongs to $r$ if and only if ones$(r) \subseteq X$ and zeros$(r) \cap X = \emptyset$, and for all $1\leq i \leq t$ it holds that
 $$\mbox{prem}[i] \subseteq X \quad \Ra \quad \ \mbox{conc}[i] \subseteq X.$$
 While this defintion reflects the author's Mathematica-implementation of $\{0,1,2,a,b\}$-valued rows in the upcoming $a,b$-algorithm, we shall use a more visual representation. Thus, up to permutation of the entries, a typical $\{0,1,2,a,b\}$-valued row looks as follows:
 
 (5) \qquad $r\quad = \quad (0,0,1,1,2,2,2,a_1, b_1, b_1, a_2, b_2, b_2, b_2, a_3, b_3)$.
 
For instance, prem$(2) = \{11\}$, conc$(2)= \{12, 13, 14\}$. It is clear  that generally
 
 (6) \qquad $|r|\quad =\quad 2^{|{\rm twos}(r)|} \cdot \ds\prod_{i=1}^t(1+2^{|{\rm prem}(i)|})$ \hspace*{2cm} $\left( \mbox{where} \ \ds\prod_{i=1}^0 (\cdots ) =1 \right)$.

Each family $\sum$ of singleton-premise implications corresponds to a directed graph $D$ in the obvious way. Let $(W, \leq )$ be the factor poset of the strong components of $D$. 
Since the $\sum$-closed sets correspond bijectively to the ideals of $(W, \leq)$, we henceforth restrict attention to posets. Coupled to a poset $(W, \leq)$ we consider the natural implicational base $\sum = \sum (W, \leq )$ that we encountered in the proof of Corollary 3. Upon relabelling we may assume that $1,2, \cdots, w$ is a linear extension of $(W, \leq)$. The benefit is that when we impose $A_{i+1} \ra B_{i+1}$ (say $\{j\} \ra B_{i+1}$) upon a multivaled row $r$, the $j$-th entry $\alpha_j$ of $r$ will always be  $\alpha_j =2$ (as opposed to Case 3 in Sec.3). Let us see in detail how $\{j\} \ra B_{i+1}$ is to be imposed on the $\{0,1,2,a,b\}$-valued row $r$. Different from Section 3 we base our case distinction on the behaviour of $B = B_{i+1}$.

{\it Case 1:} $B \cap \ \mbox{zeros}(r) = \emptyset$. Then the 2 at the $j$th position becomes $0$.

{\it Case 2:} $B \subseteq \ \mbox{ones}(r)$. Then $r$ carries over unaltered.

{\it Case 3:} At least one $B$-position is 2 and the others, if any, are 1. Then the 2 at the $j$th position becomes the symbol $a$ and the $2$'s indexed by $B$ become symbols $b$. (Both $a$ and the $b$'s get an appropriate index to distinguish them from other such symbols potentially present in $r$.)

{\it Case 4:} There are no $0$'s but symbols $a_i$ or $b_i$ within the range of $B$.  A typical situation would be:
$$\begin{array}{llcccccccccccc}
r& = & (\cdots, & b_1, & b_2, & 0, & \underbrace{b_4, b_2, a_1, a_2, b_3, 2}_{B}, & b_3, & a_3,  & a_4,  & \underbrace{{\bf 2}}_{j}, & 2, & \cdots, & 2)\\
\\
r'& = & (\cdots, & b_1, & b_2, & 0, & b_4, b_2, a_1, a_2, b_3, 2, & b_3, & a_3, & a_4, & {\bf 0}, & 2, & \cdots, & 2) \\
\\
r'' & =& ( \cdots, & 1, & 1, & 0, & 1, \ 1,  \ 1, \ 1, \ 1, \ 1, & b_3, & a_3, & 2, & {\bf 1}, &  2, & \cdots, & 2) \end{array}$$
One checks that the sets $X \in r$ that satisfy $\{j \} \ra B$ are precisely the ones in the disjoint union of $r'$ and $r''$. As to $r''$, notice that putting one $b_3 =1$ in $a_3b_3b_3$ results in $a_3b_3$, whereas putting $b_4=1$ in $a_4b_4$ results in $2$. That concludes the formal definition of the $(a,b)$-algorithm.


Lemma 1 in Section 3 was formulated for $k$-element models but mutatis mutandis also holds for models of cardinality $\leq k$ (see [W2]). In particular, if $k=w$, it is about enumerating {\it all} models. This is the version we need for assessing the $(a,b)$-algorithm. Namely, due to the linear extension ordering, all occuring rows are automatically feasible, and so $f(h,w) =0$. Furthermore $h =w$ and $s=2$. Hence the $O(Rhs(w+f(h,w)))$ term in Lemma 1 becomes $O(Rw(w+0))=O(Rw^2)$. To summarize:

\begin{tabular}{|l|} \hline \\
{\bf Theorem 5:} The $(a,b)$-algorithm takes time $O(Rw^2)$ to compactly display the ideals \\
of a $w$-element poset as a disjoint union of $R$ many $\{0,1,2,a,b\}$-valued rows.\\ \\ \hline \end{tabular}

As in all versions of POE, in theory $R$ can only be bound by the total number $N$ of models, yet in practise (as in Table 4) often $R \ll N$. From the introduction recall the $O(Nd)$ algorithm of [MN]. While $d$ is better than $w^2$, this is outweighed by $R \ll N$.

Enumerating the ideals of a poset has applications in scheduling theory and other parts of operations research [CLM].  As to applications in pure algebra, one concerns the calculation of all submodules of a finite $R$-module from its join-irreducible submodules. In fact that's only one instance of calculating a finite modular lattice from its poset of join irreducibles $p$, and a knowledge of its collinear triplets $(p_1, p_2, p_3)$ in the sense that $p_1 \vee p_2 = p_1 \vee p_3 = p_2 \vee p_3$. Here one benefits from the flexibility of the $(a,b)$-algorithm to target only specific, further constrained ideals. Albeit [SW] needs polishing, it provides the details. It is long known [GI] that all ``stable marriages'' w.r. to given rank ordered preference lists constitute a distributive lattice, and whence could be viewed as ideals of a poset and enumerated by the $(a,b)$-algorithm. Due to the flexibility of the $(a,b)$-algorithm that approach is tempting when only particular stable marriages are targeted, which do not yield to other methods. Another application [W3] is the production of all pairwise nonisomorphic incidence algebras that share some common pattern.

\section{Calculating Whitney numbers of distributive lattices}

The number $N_k$ of rank $k$ elements of any ranked lattice $L$ is called the $k$-th {\it Whitney number} of $L$. Here we focus on the distributive lattice $D(P)$ of all ideals $X$ of a $w$-element poset $W = (W, \leq )$, in which case the rank of $X$ equals its  cardinality. Thus $N_k$ can be calculated in time $O(N_kw^3)$ (Corollary 3). We refer to the introduction for various special cases concerning $N_k$. In practise the $(a,b)$-algorithm beats the $0,1,2$-algorithm underlying Corollary 3. Even more so when the issue is not enumeration but merely the calculation of $N_k$. But in order to do so the $(a,b)$-algorithm needs a little update. Namely, given a $\{0,1,2,a,b\}$-valued row $r$, we wish to efficiently calculate
 
 (7) \qquad $\mbox{Card}(r,k) \quad : = \quad  |\{X \in r: \ |X| = k\}|$
 
 for all $k \in [w]$. This is more subtle than (6). Associate with $r$ a polynomial pol$(r,x)$ with indeterminate $x$ that comes as a product as follows. Each symbol 1 contributes a factor $x$, and each symbol 2 a factor $1+x$. Further, each pattern $a_i b_i b_i \cdots b_i$ with $m = |\mbox{prem}(i)|$ contributes a factor 
 $$1+{m \choose 1}x + {m \choose 2} x^2 + \cdots + {m \choose m} x^m+x^{m+1}.$$
 It is easy to verify that Card$(r,k)$ is the coefficient of $x^k$ in pol$(r,x)$. For instance,  $r$ in (5) has
 $$\begin{array}{lll} \mbox{pol}(r,x) & = & x^2(1+x)^3(1+2x+x^2+x^3)(1+3x+3x^2+x^3+x^4)(1+x+x^2)\\
 \\
 & =& x^2+9x^3+37x^4+93x^5+\cdots + 6x^{13} +x^{14} \end{array}$$
 and thus e.g. Card$(r, 5) =93$.
   The multiplying out of factors is achieved ``fast'' using the Mathematica command {\tt Expand} and we do not to delve into its theoretical complexity here (for some $e$-algorithm handling hypergraph transversals this was done elsewhere). Summing up the polynomials of all final rows yields the {\it rank polynomial}\footnote{This is not to be confused with the {\it order polynomial} of $(W, \leq)$ whose coefficients count the number of homomorphisms from $W$ to chains.}
  $$RP(x,D)\quad  = \quad \ds\sum_{k=0}^wN_kx^k$$
  of our distributive lattice $D=D(W)$.
  
  Let us present another, now recursive way, to calculate $RP(x) = RP(x,D)$. For any $a \in W$ let $a \downarrow = \{p \in W: \ p \leq a \}$ and $a \uparrow = \{p \in W: \ p \geq a\}$ be the generated ideal respectively filter. Put $T_- = \{X \in D: \ a \not\in X\}$ and $T_+ = \{X \in D: \ a \in X\}$. One verifies that $T_-$ equals $D(W\setminus a\uparrow)$, and $X \mapsto X \setminus (a \downarrow )$ is a bijection from $T_+$ onto $D (W \setminus a \downarrow)$.  It follows\footnote{The recursion (8) has likely been discovered before but the author couldn't pinpoint a reference. We mention that a similar recurrence used in [BK, p.106] to count the antichains of $(W,\leq)$, in fact is the standard recurrence to calculate the independence polynomial of a graph (in this case the comparability graph of $(W, \leq)$. Notice that while antichains and order ideals of $(W, \leq)$ are in bijection, this is not the case for the respective $k$-element objects.} that
  
  (8) \qquad $RP(x) \quad = \quad  RP_-(x)\quad + \quad x^\alpha RP_+(x)$
  
  where $RP_-(x)$ and $RP_+(x)$ are the rank polynomials of $D(W \setminus a\uparrow)$ and $D(W \setminus a \downarrow)$ respectively, and $\alpha = |a\downarrow |$.
  
  Iterating (8) on the smaller posets $W \setminus a\downarrow$ and $W \setminus a\uparrow$ eventually yields $RP(x)$ explicitely. Here it pays to always pick $a$ in such a way that $|a \downarrow | +|a \uparrow|$ is maximum. Furthermore, for $n$-element antichains we can dispense with (8) since their rank polynomials are obvious. (That's more generally true for disjoint unions of chains but they are hard to identify.)
  
  In Table 4 below we pit the $(a,b)$-algorithm against the recursive method on some random examples. Specifically the lattices $D$ are the ideal lattices of posets $W(m, \ell, t)$  which consist of $\ell$ levels Lev$(i)$, all of cardinality $m$, in such a way that each $a \in \, \mbox{Lev}(i) \ (2 \leq i \leq \ell)$ is assigned $t$ random lower covers in Lev$(i-1)$.
  \hfill $\blacksquare$

Both algorithms were implemented with Mathematica 8.0, the implementation of the recursive method being straightforward\footnote{It would be harder to set up and implement recursions for counting the $k$-element subtrees of a tree; cf. Corollary 4.}. In all cases both methods agreed upon $RP(x,D)$, thus very likely both are correct. Rather than $RP(x,D)$ we record the total number $N = RP(1,D)$ of ideals, the respective times, as well as the parameters $R$ and $nsum$. The former is the number of final rows generated by the $(a,b)$-algorithm (matching the notation of Theorem 5), the latter the number of final (i.e. coupled to antichains) summands constituting $RP(x,D)$ when using the recursive method. Depending on the type of poset one or the other algorithm prevails. Recall that the $(a,b)$-algorithm has the additional benefit of delivering the ideals {\it themselves} in a compact format.

\begin{tabular}{c|c|c|c|c|c} 
$(m, \ell, s)$ & $N$ & $R$ & time-$(a,b)$ &  $nsum$ & time-rec \\ \hline
$(5,28,2)$ & $30416508$ &  $8932$ & 25 & $4827$ & $7$\\ \hline 
$(3,28,1)$ & $66506610198$ & $56256$ & 61 & $595830$ & $1070$\\ \hline
$(10, 15, 9)$ & $15560$ &  $984$ & 3 & $228$ & $1$\\ \hline
$(15, 10, 6)$ & $709150$ &  $25302$ & 80 & $5870$ & $6$\\ \hline
$(20,2,2)$ & $500466940$ & $1872$ & 3 & $30205$ & $49$ \\ \hline
$(20,3,2)$ & $27003561664$ & $1041866$ & 746 & $946439$ & $1612$\\ \hline
$(20,4,2)$ & $214942410560$ & $1148952$ & 1335 & $6713622$ & $11892$ \\ \hline
\end{tabular}

\hspace*{2cm} Table 4

\section{About quasi-implications}

Each singleton-premise implication, say $\{a\} \ra \{b_1, b_2, b_3\}$ properly speaking is the propositional formula in (i) below. It serves as point of departure in our attempt to carry over the $a,b$-algorithm to other kinds of (singleton-premise) {\it quasi-implications}. We investigate the eight variations obtained from (i) by independently negating the premise $a$, respectively negating all variables $b_1, b_2, b_3$ of the conclusion, respectively switching $\wedge$ to $\vee$. Since negating {\it all} variables can't deliver anything logically new (yet psycho-logically), we mainly focus on (i) to (iv) rather than (i$^\ast$) to (iv$^\ast$):
\begin{enumerate}
	\item [(i)] \ $a \ra (b_1 \wedge b_1 \wedge b_3)$  \hspace*{2cm} $(1,1)$ implication
	\item[(ii)] \ $\ol{a} \ra (\ol{b}_1 \vee \ol{b}_2 \vee \ol{b}_3)$ \hspace*{2cm} $(0,0')$ quasi-implication
	\item[(iii)] \ $a \ra (\ol{b}_1 \vee \ol{b}_2 \vee \ol{b}_2)$ \hspace*{2cm}  $(1,0')$ quasi-implication
	\item[(iv)]  \ $a \ra (\ol{b}_1 \wedge \ol{b}_2 \wedge \ol{b}_3)$ \hspace*{2cm} 	$(1,0)$ quasi-implication
	
	\item [(i$^\ast)$] \ $\ol{a} \ra (\ol{b}_1 \wedge \ol{b}_2 \wedge \ol{b}_3)$ \hspace*{2cm} $(0,0)$ quasi-implication
	\item[(ii$^\ast)$] \ $a \ra (b_1 \vee b_2 \vee b_3)$ \hspace*{2cm} $(1,1')$ quasi-implication
	\item[(iii$^\ast)$] \ $\ol{a} \ra (b_1 \vee b_2 \vee b_3)$ \hspace*{2cm} $(0,1')$ quasi-implication
	\item[(iv$^\ast)$]  \ $\ol{a} \ra (b_1 \wedge b_2 \wedge b_3)$ \hspace*{2cm} $(0,1)$ quasi-implication
\end{enumerate}

As to (ii), this is equivalent to $(b_1 \wedge b_2 \wedge b_3)\ra a$. Such singleton-{\it conclusion} implications are more powerful than singleton-premise implications since {\it each} implication, say $(b_1 \wedge b_2 \wedge b_3) \ra (a_1 \wedge a_2)$, is equivalent to a conjunction of such implications $(b_1 \wedge b_2 \wedge b_3) \ra a_1$ and $(b_1 \wedge b_2 \wedge b_3)\ra a_2$. 
The $(a,b)$-symbolism does not\footnote{The author's attempt to carry it over degenerated in a messy case distinction and thus spawned the implication $n$-algorithm touched upon at the beginning of Section 5.} carry over to implications with premises of cardinality $\geq 1$.

As to (iii), this is equivalent to the clause $\ol{a} \vee \ol{b}_1 \vee \ol{b}_2 \vee \ol{b}_3$. Also (ii) could have been written as $a\vee \ol{b}_1 \vee \ol{b}_2 \vee \ol{b}_3$. By definition a conjunction of clauses with at most one positive literal, is called a {\it Horn formula}. The {\it Horn $n$-algorithm} of [W2] extends the implication $n$-algorithm in obvious ways. Having at most one {\it negative} literal per clause yields {\it Anti-Horn formulae} (cases (ii$^\ast)$ and (iii$^\ast)$).

As to (iv),  consider the task to determine all independent sets of vertices (anticliques) in a graph $G$. Obviously $X$ being an anticlique amounts to the satisfisfaction of all clauses $\ol{a} \vee \ol{b}$ where $\{a,b\}$ ranges over the edges of $G$. Grouping together all edges incident with a fixed vertex $a$ this becomes (say)
$$(\ol{a} \vee \ol{b}_1) \wedge (\ol{a} \vee \ol{b}_2) \wedge (\ol{a} \vee \ol{b}_3) \quad \equiv \quad \ol{a} \vee (\ol{b}_1 \wedge \ol{b}_2 \wedge \ol{b}_3)\quad \equiv \quad a \ra (\ol{b}_1 \wedge \ol{b}_2 \wedge \ol{b}_3)$$
This gives rise to the $(a, \ol{b})${\it -algorithm} (work in progress) which bears a superficial resemblence to the $(a,b)$-algorithm. However, (i) and (iv) are not duals in any sense, and also the technical details differ considerably. For instance, there is no such thing as a linear extension for graphs.  

Apart from (i) to (iv) versus (i$^\ast$) to (iv$^\ast$), one can distinguish the quasi-implications that feature $\wedge$ in their conclusions of fixed type (carrying tags $(x,y)$) from the ones that feature $\vee$ (carrying tags $(x,y')$). Acyclic families of quasi-implications invite conglomeration to rooted trees, and this has indeed algorithmic benefits for type $(x,y)$. Let us sketch the basic idea for type $(1,1)$. The conjunction

(9) \qquad $(1 \ra (2 \wedge 3 \wedge 4)) \ \wedge \ (2 \ra (5 \wedge 6)) \ \wedge \ (4 \ra 7) \ \wedge \ (5 \ra (8 \wedge 9))$

of four implications gives rise to the (rooted) tree in Figure 2(a).
 
\begin{center}
\includegraphics{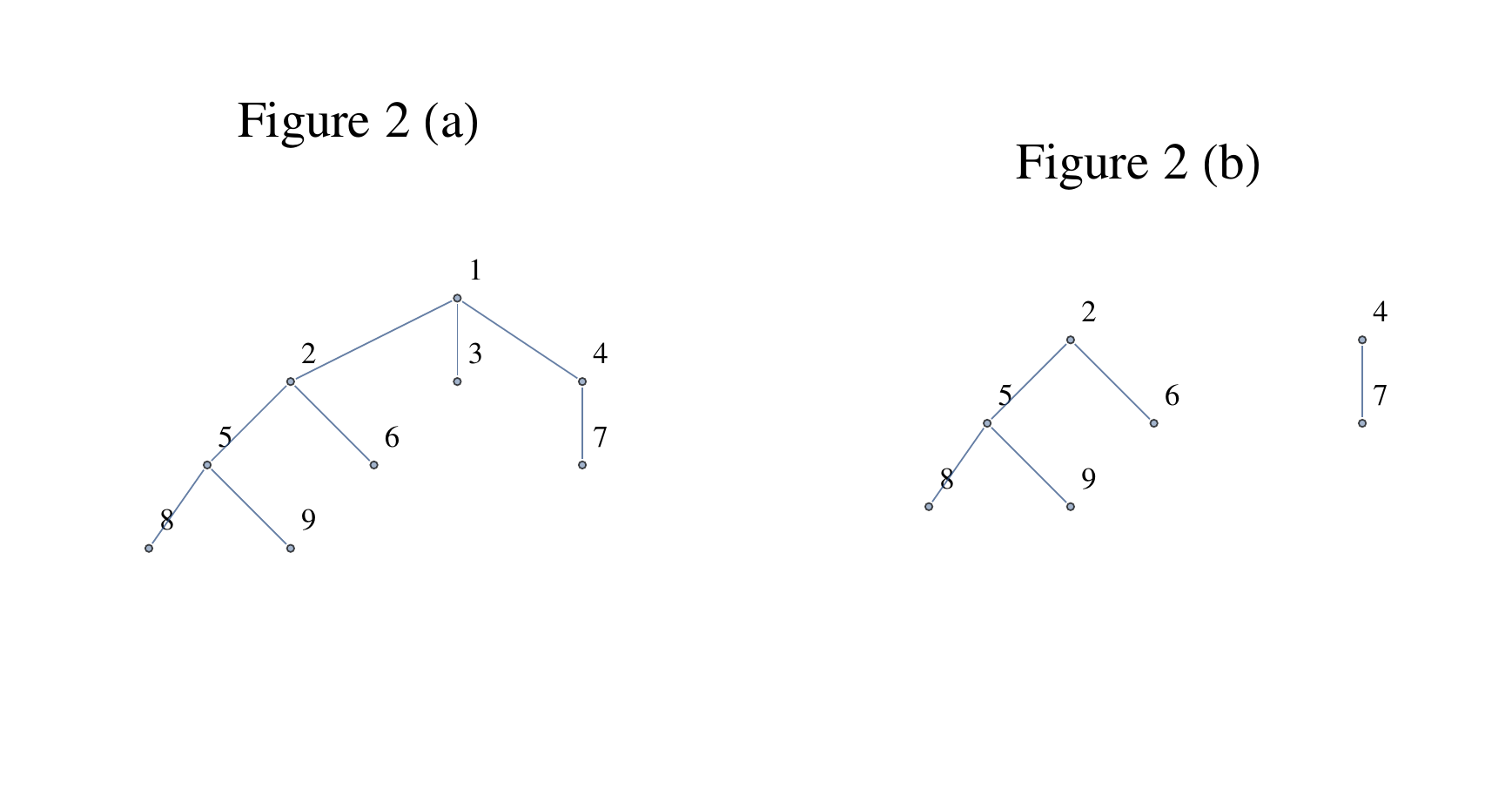}
\end{center}

We wish to encode with appropriate symbols $\ast$ the set of all $X \in {\cal P}([9])$ that satisfy (10) as just {\it one} ``tree-valued'' row $r = (\ast, \ast, \ast, \ast, \ast, \ast, \ast, \ast, \ast )$. Notice that $\{X \in r: \ 1 \in X\} = (1,1,1,1,1,1,1,1,1)$ whereas $\{X \in r: \ 1 \not\in X\} = (0, \ol{\ast}, 2, a, \ol{\ast}, \ol{\ast}, b, \ol{\ast}, \ol{\ast})$ (see Figure 2(b)). Details pending, the aim is to impose trees $T$ upon tree-valued rows $r'$ by splitting $r'$ in $r'_1$ and $r'_2$ according to whether $T$'s root is turned to $1$ or $0$. The resulting ``tree-algorithm'' would comprise our $(a,b)$-algorithm as the special case where all occuring trees have height 1. 

We indicated algorithms for handling all eight types individually, but what if some types occur simultaneously? A case in point is [W1] where (i), (i$^\ast$), (iv), (iv$^\ast$) do just that.  Here an ad hoc algorithm of the $0,1,2$-type is successful.

\section*{References}
\begin{enumerate}
\item[{[BK]}] J. Berman, P. K\"{o}hler, Cardinalities of finite distributive lattices, Mittelungen Math. Sem. Giessen, Heft 121 (1976) 103-124.
\item[{[C]}] A. Conflitti, On Whitney numbers of the order ideals of generalized fences and crowns, Disc. Math. 309 (2009) 615-621.

\item[{[CLM]}] N. Caspard, B. Leclerc, B. Monjardet, Finite ordered sets, Encyclopedia Math. and Appl. 144, Cambridge University Press 2012.
\item[{[[EJ]}] P.H. Edelman, R.E. Jamison, The theory of convex geometries, Geom. Dedicata 19 (1985) 247-270.
\item[{[FK]}] U. Faigle, W. Kern, Computational  complexity of some maximum average weight problems with precedence constraints, Oper. Res. 42 (1994) 688-693.
\item[{[GI]}] D. Gusfield, R.W. Irwing, The stable marriage problem: Structure and algorithms, MIT Press 1989.
		\item[{[HNS]}] M. Habib, L. Nourine, G. Steiner, Gray codes for the ideals of interval orders, Journal of Algorithms 25 (1997) 52-66.
	\item[{[K]}] D. Knuth, The Art of Computer Programming, Vol.4, Fasc.3, 2009, Pearson Educ. Inc. 
\item[{[MN]}] R. Medina, L. Nourine, Algorithme efficace de g\'{e}n\'{e}ration des ideaux d'un ensemble ordonn\'{e}, C.R. Acad. Sci. Paris S\'{e}r. I Math. 319 (1994) 1115-1120.		
\item[{[R]}] F. Ruskey, Listing and counting subtrees of a tree, SIAM J. Comput. 10 (1981) 141-150.
\item[{[RM]}] H. Mannila, K.J. R\"{a}ih\"{a}, The design of relational databases, Addison-Wesley 1992.

		\item[{[S]}] G. Steiner, On estimating the number of order ideals in partial orders, with some applications, J. Stat. Planning and Inference 34 (1993) 281-290.
		\item[{[SW]}] Y. Semegni, M. Wild, Lattices freely generated by posets within a variety. Part II: Finitely generated varieties, arXiv:1007.1643.
			\item[{[W1]}] M. Wild, Revisiting the enumeration of all models of a Boolean 2-CNF, arXiv: 1208.2559
	\item[{[W2]}] M. Wild, Compactly generating all satisfying truth assignments of a Horn formula, Journal of Satisfiability, Boolean  Modeling and Computation, 8 (2012) 63-82.
	\item[{[W3]}] M. Wild, Incidence algebras that are uniquely determined by their zero-nonzero matrix pattern, Lin. Alg. and Appl. 430 (2009) 1007-1016.
\end{enumerate}

\end{document}